\documentclass[journal ]{new-aiaa}
\usepackage[utf8]{inputenc}
\usepackage{textcomp}

\usepackage{graphicx}
\usepackage{amsmath}
\usepackage{mathtools}
\usepackage[version=4]{mhchem}
\usepackage{siunitx}
\usepackage{longtable,tabularx}
\newcommand{\hide}[1]{}

\usepackage{xcolor}
\definecolor{cardinal}{rgb}{0.0, 0.8, 0.6}

\newcommand{\bx}{\textbf{\emph{x}}}

\newcommand{\ba}{\textbf{\emph{a}}}
\newcommand{\bb}{\textbf{\emph{b}}}
\newcommand{\bc}{\textbf{\emph{c}}}
\newcommand{\bd}{\textbf{\emph{d}}}
\newcommand{\bl}{\boldsymbol{\ell}}

\newcommand{\bw}{\textbf{\emph{w}}}

\newcommand{\br}{\textbf{\emph{r}}}
\newcommand{\bs}{\textbf{\emph{s}}}

\newcommand{\bv}{\textbf{\emph{v}}}

\newcommand{\bp}{\textbf{\emph{p}}}
\newcommand{\bq}{\textbf{\emph{q}}}

\newcommand{\bT}{\textbf{\emph{T}}}

\newcommand{\PP}{{\mathbb P}}
\newcommand{\RR}{{\mathbb R}}

\title{Initial Orbit Determination from Only Heading Measurements}

\author{John A. Christian \footnote{Associate Professor, Guggenheim School of Aerospace Engineering. Associate Fellow AIAA.}}
\affil{Georgia Institute of Technology, Atlanta, GA 30332}

\begin{document}

\maketitle

\begin{abstract}
This work introduces the problem of initial orbit determination (IOD) from only heading measurements. Such a problem occurs in practice when estimating the orbit of a spacecraft using visual odometry measurements from an optical camera. After reviewing the problem geometry, a simple solution is developed in the form of an iterative scheme on the parameters describing the orbital hodograph. Numerical results are presented for an example spacecraft in low lunar orbit. The principal intent of this brief study is to communicate the existence of a new class of IOD problem to the community and to encourage the broader study of hodographs and heading-only IOD.
\end{abstract}

\section{Introduction}

Orbit determination from a single measurement type is one of the classical problems of astrodynamics. Before the era of space exploration, observations of celestial bodies were exclusively collected by ground-based telescopes. When using such optical instruments for orbit determination, the observations usually take the form of a direction measurement that form a line connecting the locations of the observer (e.g., telescope) and the orbiting body. Thus, if the orbiting body and observer have positions $\br_i$ and $\bp_i$ at time $t_i$, then the observed direction $\bl_i$ is given by
\begin{equation}
    \bl_i \propto \br_i - \bp_i
\end{equation}
where $\bl_i \in \PP^2$ and has ambiguous scale. For more details on projective space (e.g., $\PP^2$, $\PP^n$) see Refs.~\cite{Hartley:2003,Gallier:2011,Christian:2021d}.
It is unsurprising, therefore, that the earliest work in IOD considered orbit reconstruction with only direction (i.e., line-of-sight, bearing, or angle) measurements. Given three generic direction measurements $\{ \bl_i \}_{i=1}^3$ at known times $\{ t_i \}_{i=1}^3$, one may solve for the orbit using the classical methods of Laplace \cite{Laplace:1780}, Gauss \cite{Gauss:1809}, Gooding \cite{Gooding:1997}, and others \cite{Escobal:1976}. In the special case of all the observations being coplanar, four observations are required \cite{Baker:1977}. 

Other classical IOD problems make use of position vectors instead of directions. Amongst the earliest of these is the so-called Lambert's problem \cite{Lambert:1761,Volk:1980}, which describes the task of finding the orbit from two position vectors $\{\br_i\}_{i=1}^2$ and time of flight between them $\Delta t_{12} = t_2 - t_1$. Lambert's problem was first practically solved by Lagrange in 1773 \cite{Lagrange:1773,Albouy:2019} and has since been widely studied \cite{Gooding:1990,Russell:2019}. The other important position-only problem is Gibbs problem \cite{Gibbs:1889}, which describes the task of finding the orbit from three position vectors $\{\br_i\}_{i=1}^3$. Gibbs problem does not require explicit knowledge of the observation times. The ultimate objective for both Lambert's problem and Gibbs problem is to recover the unknown velocity vector corresponding to one of the known position vectors.

The principal narrative within the contemporary IOD literature still focuses on observations of either direction or position. Although some IOD methods now include other measurements (e.g., angle rates, range, range-rate), the fundamental problems remain largely unchanged from the classical results mentioned above. In modern practice, IOD methods usually provide a preliminary orbit solution for one of two purposes. The first purpose is acquisition of an object's orbit with minimal observations so that sensors may be tasked to collect additional observations (or so that existing unlabeled observations may be properly attributed to this object). The second purpose of IOD methods is to provide an initial guess to a precise orbit determination (POD) method. Such POD methods usually processes a large number of measurements and consider non-Keplerian motion (e.g., non-spherical gravity, atmospheric drag, solar radiation pressure). Moreover, most POD methods produce a maximum likelihood estimate of the orbit \cite{Tapley:2004}, which is not generally the case for IOD solutions. 

Until recently, the problems described above (and their close relatives) have constituted the extent of the mainstream study of IOD. However, it was suggested in 2018--2019 that one might reimagine the problems of Lambert and Gibbs to be given velocity vectors instead of position vectors \cite{Christian:2018,Hollenberg:2019}. This gives rise to the problem of IOD from two velocity vectors and time-of-flight (analog to classical Lambert's problem) and of IOD from three velocity vectors (analog to classical Gibbs problem). Efficient geometric solutions are possible using the orbital hodograph. These velocity-only IOD problems arise with sensors capable of directly measuring velocity, such as X-Ray pulsar navigation (XNAV) \cite{Hou:2022} or StarNAV \cite{Christian:2019}. The ultimate objective is to recover the unknown position vector corresponding to one of the known velocity vectors.

If the work of Hollenberg and Christian from Refs.~\cite{Christian:2018} and \cite{Hollenberg:2019} introduced the problem of IOD from velocity observations, it follows that one might also consider IOD from velocity directions (i.e., from headings). In this work, the \emph{heading} of an orbiting body at time $t_i$ is denoted as $\bs_i \in \PP^2$ and describes the direction of the inertial velocity $\bv_i = \dot{\br}_i$. Thus,
\begin{equation}
    \bs_i \propto \bv_i \propto \bv_i / \| \bv_i \|
\end{equation}
where $\bs_i$ is a $3 \times 1$ vector of arbitrary scale. Choosing $\bs_i$ to be a unit vector (i.e., to be $\bv_i / \| \bv_i \|$) is often convenient, but not required for the algorithms developed here. IOD with only heading measurements is the analog of IOD with only direction measurements (e.g., with classical solutions of Laplace, Gauss, Gooding, and so on)---in this new problem, however, one observes the direction along the velocity vector instead of the direction along the relative position. A solution from only inertial headings $\bs_i$ requires four measurements since the measurements are all coplanar (for the same reasons that four measurements are needed with coplanar directions $\bl_i$, see Ref.~\cite{Baker:1977}). Thus, this work proposes a new class of IOD problem:\\ \\

\hangindent=0.5cm
\textbf{Problem Statement:} \emph{Consider an object $\mathcal{A}$ following an unknown Keplerian orbit about a central body $\mathcal{B}$ of known mass. Determine the orbit of $\mathcal{A}$ relative to $\mathcal{B}$ given four observations of the object's heading $\{ \bs_i \}_{i=1}^4$ collected at four known times $\{ t_i \}_{i=1}^4$.}
\\

The problem of IOD from headings is not a purely academic question, but arises in a number of practical scenarios. As an example, consider a spacecraft orbiting a body (e.g., Earth, Moon) equipped with a camera for navigation. Suppose this camera is pointed towards the central body and captures images of the surface. If surface features in sequential images are matched to one another, visual odometry may be used to infer the relative motion between the camera and surface to an unknown scale. This is advantageous since tracking features from image-to-image does not require a model (i.e., map) of the observed body. Since attitude is usually known from other sources (e.g., from a star tracker aboard the spacecraft), the visual odometry algorithm provides a surface-relative heading measurement. If the orbital velocity is much larger than the velocity of the surface points, then the surface-relative heading is usually very close to $\bs_i$. Such heading measurements may be efficiently computed with the algorithms of Ref.~\cite{Christian:2021b}, which refers to $\bs_i$ as the ``direction-of-motion'' (the same as heading when the time between an image pair is small). Therefore, solutions to the heading-only IOD problem provide a practical means for map-free orbit determination with only a single camera, with the measurements coming from a visual odometry algorithm.

The sections that follow consider the geometry and dynamics of heading measurements, which are easiest to understand with the orbital hodograph. These geometric insights are used to create a numerical algorithm to the heading-only IOD problem and a simple experiment is presented to evaluate performance. The development of more nuanced solution schemes is the topic of ongoing and future work.

\section{The Orbital Hodograph}
Consider an object orbiting a large central body. Define the \emph{orbital hodograph} as the locus of points traced out by the tip of the orbiting object's velocity vector while keeping the tail fixed at the origin. It was shown by Hamilton in 1847 \cite{Hamilton:1847} that the hodograph is always a perfect circle for a Keplerian (two-body) orbit. That is, the hodograph remains a circle regardless of the conic describing the orbit itself (e.g., circle, ellipse, parabola, hyperbola). The simple geometry of a circle makes the hodograph a useful tool for orbit analysis, as notably explored by Eades \cite{Eades:1968} and Altaman \cite{Altman:1965,Altman:1965b} in the 1960s.  More recently, the hodograph is found to be an especially powerful tool for IOD with velocity (e.g.,  Refs.~\cite{Hollenberg:2019} and \cite{Hou:2022}) and heading (e.g., this work) observations. A brief development of the orbital hodograph is now presented.

Suppose $\br \in \RR^3$ is the position of a spacecraft relative to a large central body. Classical two-body orbital mechanics leads to the equations of motion 
\begin{equation}
    \ddot{\br} = - \frac{\mu}{r^3} \br
\end{equation}
where $\mu$ is the central body's gravitational parameter and $r = \| \br \|$. The solution to this differential equation within the orbit plane at time $t_i$ may be written in polar coordinates as
\begin{equation}
    r_i = \frac{h^2 / \mu}{1 + e \cos \theta_i}
\end{equation}
where $h$ is the specific angular momentum, $e$ is the orbit eccentricity, and $\theta_i$ is the true anomaly.

Let $\{ \hat{\bp},\hat{\bq},\hat{\bw} \}$ describe the right-handed orthonormal basis defining the perifocal frame. Choose $\hat{\bp}$ to be the unit vector from the central body towards orbit periapsis, $\hat{\bw}$ to be in the direction of the specific angular momentum (i.e., normal to the orbit plane), and $\hat{\bq}$ to complete the right-handed system. Using the perifocal basis, one can directly write the three-dimensional (3-D) position vector as
\begin{equation}
    \br_i = r_i \cos \theta_i \, \hat{\bp} + r_i \sin \theta_i \, \hat{\bq}
\end{equation}
Taking the time derivative of the position yields
\begin{equation}
    \label{eq:VelPerifocal01}
    \bv_i = \dot{\br}_i = \frac{\mu}{h} \left[  \left(-\sin \theta_i \right)  \hat{\bp} + \left(e + \cos \theta_i \right) \hat{\bq} \right]
\end{equation}
which describes a circle with radius $\mu/h$ and center at $\mu e / h \, \bq$. The relation in Eq.~\eqref{eq:VelPerifocal01} is a standard result found in most contemporary astrodynamics texts \cite{Wiesel:1997,Vallado:2007,Curtis:2020,Bate:2020}. Interestingly, few texts mention that Eq.~\eqref{eq:VelPerifocal01} implies that the velocity follows a circular path and that this circle is called the orbital hodograph.

\begin{figure}[t!]
    	\centering
    	\includegraphics[width=0.4 \textwidth]{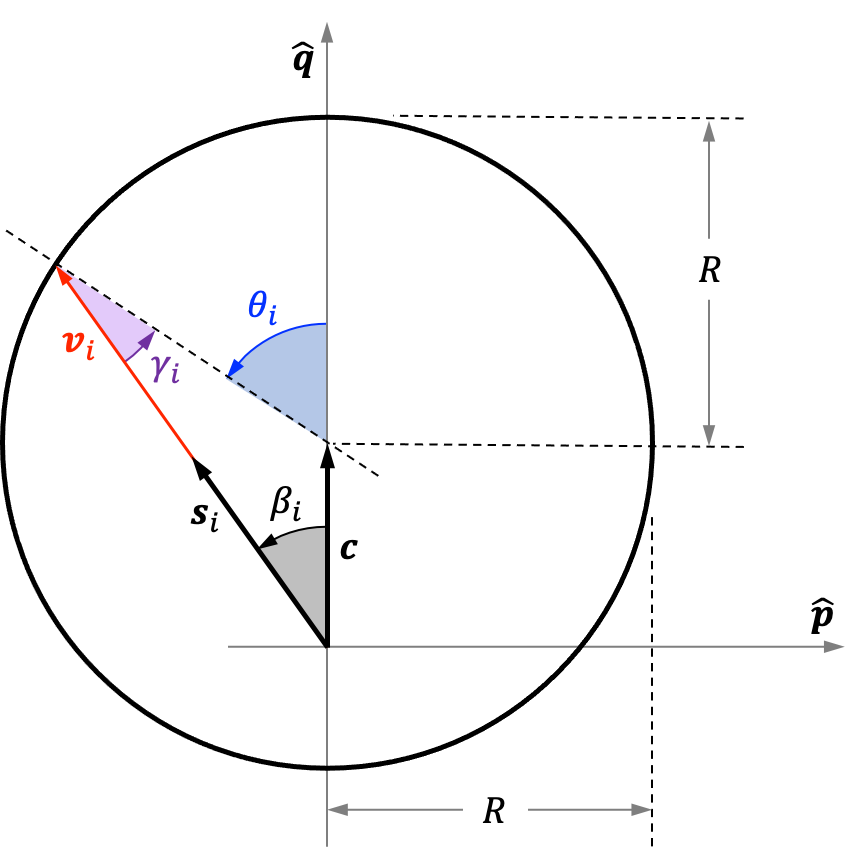}
    	\caption{Illustration of the hodograph.}
    	\label{fig:hodograph}
\end{figure}

Therefore, define the hodograph parameters $R$ and $\bc$ (which describe the circle radius and center coordinates) to be
\begin{subequations}
\label{eq:HodographParam}
\begin{equation}
    \label{eq:HodographRadius}
    R = \mu / h
\end{equation}
\begin{equation}
    \label{eq:HodographCenter}
    \bc = (\mu e / h ) \, \hat{\bq} = R e \hat{\bq}
\end{equation}
\end{subequations}
such that the velocity may be written as
\begin{equation}
    \label{eq:VelVecHodograph}
    \bv_i = R \left[  \left(-\sin \theta_i \right)  \hat{\bp} + \cos \theta_i  \hat{\bq} \right] + \bc
\end{equation}
This results in the hodograph geometry illustrated in Fig.~\ref{fig:hodograph}. Also labeled in this figure is the flight path angle, $\gamma_i$, which is defined as the angle between the velocity and the local horizontal. The geometry in Fig.~\ref{fig:hodograph} may also be used to introduce the \emph{heading angle}, $\beta_i$, defined as the angle between the $\hat{\bq}$ direction and the velocity vector (positive in the direction of orbital motion). Thus, using Fig.~\ref{fig:hodograph} and substituting from Eq.~\eqref{eq:VelPerifocal01}, the heading angle may be computed as
\begin{equation}
    \label{eq:DefHeadingAngle}
    \tan \beta_i = \frac{-\hat{\bp}^T \bs_i}{ \hat{\bq}^T \bs_i} = \frac{-\hat{\bp}^T \bv_i}{ \hat{\bq}^T \bv_i} = \frac{\sin \theta_i}{e + \cos \theta_i}
\end{equation}
which will soon become an important relationship.

\subsection{Classical Orbital Elements from Hodograph Parameters}
There is a one-to-one mapping between the hodograph parameters and the in-plane orbital elements. To see this, observe that the results of Eq.~\eqref{eq:HodographParam} may be rearranged to find the specific angular momentum $h$ and the eccentricity $e$
\begin{equation}
    h = \mu / R
\end{equation}
\begin{equation}
    \label{eq:EccFromHodograph}
    e = \| \bc \| / R
\end{equation}
The expressions for $h$ and $e$, may then be used to find the orbit's semi-major axis $a$
\begin{equation}
    \label{eq:SemiMajorAxisHodograph}
    a = \frac{h^2 / \mu}{1-e^2} = \frac{\mu/R^2}{1-\bc^T \bc/R^2} = \frac{\mu}{R^2-\bc^T \bc}
\end{equation}
The orbital elements $a$ and $e$ fully describe the size and shape of the orbit within the perifocal frame.

\subsection{Kepler's Equation and the Hodograph}
The IOD problem considered in this work requires a means of predicting how an orbiting body's heading changes with time. This may be achieved by manipulation of Kepler's equation.
Recall that Kepler's equation is given by
\begin{equation}
    \label{eq:KeplersEqn}
    M_i = n(t_i - t_0) = E_i + e \sin E_i
\end{equation}
where $E_i$ is the eccentric anomaly, $M_i$ is the mean anomaly, $n$ is the mean motion, and $t_0$ is the time of periapsis passage. The immediate objective is to express $n$ and $E_i$ entirely in terms of the hodograph parameters ($R$ and $\bc$) and the heading angle $\beta_i$.
The desired equation for mean motion may be found by direct substitution of Eq.~\eqref{eq:SemiMajorAxisHodograph} into the standard definition of $n$
\begin{equation}
    \label{eq:MeanMotion}
    n = \sqrt{\mu / a^3 } = (R^2-\bc^T \bc)^{3/2} / \mu
\end{equation}
which is consistent with Ref.~\cite{Christian:2021c}

More interesting is the relation between heading angle and eccentric anomaly. Recall that the usual equations for eccentric anomaly $E_i$ in terms of the true anomaly $\theta_i$ are given by \cite{Curtis:2020}
\begin{subequations}
\begin{equation}
    \cos E_i = \frac{e + \cos \theta_i}{1 + e \cos \theta_i}
\end{equation}
\begin{equation}
    \sin E_i = \frac{\sqrt{1-e^2}\sin \theta_i}{1 + e \cos \theta_i}
\end{equation}
\end{subequations}
such that
\begin{equation}
    \tan E_i = \frac{\sin E_i}{\cos E_i} = \frac{\sqrt{1-e^2} \sin \theta_i}{e + \cos \theta_i}
\end{equation}
Substitution of Eq.~\eqref{eq:DefHeadingAngle} into the expression for $\tan E_i$ yields a surprisingly simple relation between $E_i$ and $\beta_i$
\begin{equation}
    \label{eq:EccAnomalyAndHeadingAngle01}
    \tan E_i = \sqrt{1-e^2} \, \tan \beta_i
\end{equation}
which is consistent with a similar result from Ref.~\cite{Eades:1968}

It is possible to take the expression for $\tan E_i$ one step further by writing Eq.~\eqref{eq:EccAnomalyAndHeadingAngle01} directly in terms of the hodograph parameters ($R$ and $\bc$), the orbit normal direction $\hat{\bw}$, and a particular heading vector $\bs_i$. Since both $\bc$ and $\bs_i$ lie in the orbital plane and are perpendicular to $\hat{\bw}$, observe from Fig.~\ref{fig:hodograph} that 
\begin{subequations}
\label{eq:BetaTrigVecs}
\begin{equation}
    \| \bc \| \, \| \bs_i \| \sin \beta_i = \hat{\bw}^T (\bc \times \bs_i)
\end{equation}
\begin{equation}
    \| \bc \| \, \| \bs_i \| \cos \beta_i = \bc^T  \bs_i
\end{equation}
\end{subequations}
which also makes use of $\| \hat{\bw} \| = 1$. Now, recalling $\tan \beta_i = \sin \beta_i / \cos \beta_i$ and substituting the results of Eq.~\eqref{eq:BetaTrigVecs} into Eq.~\eqref{eq:EccAnomalyAndHeadingAngle01},
\begin{align}
    \tan E_i = \sqrt{1-e^2} \,  \frac{\hat{\bw}^T ( \bc \times \bs_i) }{\bc^T \bs_i} 
\end{align}
where it is observed that any non-zero scaling of the heading vector $\bs_i$ provides the same estimate of eccentric anomaly (i.e., both $\bs_i$ and $\alpha_i \bs_i$ produce the same value for $E_i$). As mentioned briefly before, selecting a scaling of $\|\bs_i\| = 1$ is convenient but not necessary.
The above result may be written entirely in terms of the hodograph parameters by substituting for the eccentricity $e$ from Eq.~\eqref{eq:EccFromHodograph},
\begin{align}
    \tan E_i = \frac{ \sqrt{R^2-\bc^T \bc} \, \left[ \hat{\bw}^T ( \bc \times \bs_i) \right] }{R \, \bc^T \bs_i}
\end{align}
Consequently, the quadrant discriminating estimate of eccentric anomaly $E_i$ is given by
\begin{equation}
    \label{eq:EccAnomalyHodographATAN2}
    E_i = \text{atan2}\left[ \left(\sqrt{R^2-\bc^T \bc} \, \left[ \hat{\bw}^T ( \bc \times \bs_i) \right] \right) ,\, \left( R \, \bc^T \bs_i \right) \right]
\end{equation}
which only depends on the hodograph parameters ($R$ and $\bc$), the orbit normal direction $\hat{\bw}$, and a heading vector $\bs_i$.

\section{Orbit Determination from Headings}
\label{Sec:IODMethod}
A simple numerical scheme may be constructed to estimate the orbit from four (or more) heading angles. 
More elegant schemes are certainly possible and the technique that follows is primarily meant to demonstrate that the heading-only IOD problem is solvable.

The procedure that follows has two primary steps. The first step is to compute the orbital plane and convert the 3-D problem into a two-dimensional (2-D) problem. The second step is to find the hodograph parameters within this plane. These two steps are now discussed in detail.

\subsection{Finding the Orbital Plane}
\label{Sec:InitialGuess}
The orbital plane may be found directly from the heading observations. Since all of the heading vectors lie within the orbital plane, two non-collinear headings span the plane of the orbit. 
Thus, choosing to define the orbit plane by its unit normal vector $\hat{\bw}$, observe that $\bs^T_i \hat{\bw} = 0$ for all four (or more) heading vectors. Stacking these equations leads to a least-squares problem of the form
\begin{equation}
   \begin{bmatrix}
        \bs^T_1 \\
        \bs^T_2 \\
        \bs^T_3 \\
        \bs^T_4 
   \end{bmatrix} \hat{\bw} = \textbf{0}_{4 \times 1}
\end{equation}
This is a null space problem and the direction $\hat{\bw}$ may be found in the least-squares sense via the singular value decomposition (SVD).

The $\hat{\bp}$ and $\hat{\bq}$ directions of the perifocal frame within the orbit plane cannot yet be determined (this comes later). For now, it is necessary to chose an intermediate basis for describing points within the orbit plane to make the subsequent computations explicit. Therefore, without loss of generality, define the right-handed orthonormal basis $\{ \hat{\ba},\hat{\bb},\hat{\bw} \}$, where
\begin{equation}
   \hat{\bb} = \frac{ \bs_i \times \hat{\bw} }{\| \bs_i \times \hat{\bw} \|} \quad \text{and} \quad \hat{\ba} = \hat{\bb} \times \hat{\bw}
\end{equation}
which will always be well-defined. Any of the heading vectors may be chosen to generate $\hat{\bb}$. Use these basis vectors to construct a rotation matrix from the inertial frame where the observations $\{ \bs_i \}_{i=1}^4$ are known to the newly constructed intermediate frame
\begin{equation}
   \bT = \begin{bmatrix}
        \hat{\ba}^T \\ \hat{\bb}^T \\ \hat{\bw}^T
   \end{bmatrix}
\end{equation}
Such that,
\begin{equation}
   \begin{bmatrix}
        s'_{i_1} \\ s'_{i_2} \\ 0
   \end{bmatrix} =
   \bs'_i = \bT \bs_i
\quad \text{and} \quad
    \begin{bmatrix}
        c'_1 \\ c'_2 \\ 0
   \end{bmatrix} =
   \bc' = \bT \bc
\end{equation}
The reverse mapping for the hodograph center $\bc$ may be efficiently computed as 
\begin{equation}
    \label{eq:CenterInverseMapping}
    \bc = \bT^T \bc' = c'_1 \hat{\ba} + c'_2 \hat{\bb}
\end{equation}

\subsection{Finding the Hodograph Parameters}
Given the basis $\{\hat{\ba},\hat{\bb},\hat{\bw}\}$ and $\{ \bs'_i \}_{i=1}^4$ from the prior section, the problem is now reduced to determining the orbit within the plane. This is accomplished with an iterative nonlinear least-squares scheme to estimate the three in-plane hodograph parameters $R$, $c'_1$, and $c'_2$. The construction of this scheme has two parts: (1) selection of a reasonable initial guess and (2) construction of the iterative scheme.

\subsubsection{Selecting an Initial Guess}
Before beginning the iterative scheme, one must first find a principled way of computing an initial guess. The simplest choice---and one that permits a direct solution---is to assume a circular orbit. If the orbit is circular, then $e=0$ and it follows from Eq.~\eqref{eq:HodographCenter} that $c'_1 = c'_2 = 0$. Thus the only task is to compute an initial guess for $R$.

If the orbit is circular then $\beta_i = \theta_i = E_i = M_i$. From this, and assuming $t_j > t_i$, it follows from Kepler's second and third laws that
\begin{equation}
    \frac{t_j - t_i}{\beta_j - \beta_i} = \frac{t_j - t_i}{\theta_j - \theta_i} = \frac{P}{2 \pi} = \sqrt{\frac{r_0^3}{\mu}}
\end{equation}
where $P$ is the orbit period and $r_0$ is the radius of the corresponding circular orbit. And, therefore,
\begin{equation}
    r^3_0 = \mu \left( \frac{t_j - t_i}{\beta_j - \beta_i} \right)^2 
\end{equation}
Since $\bc = \textbf{0}_{3\times 1}$ for a circular orbit, Eq.~\eqref{eq:SemiMajorAxisHodograph} shows that $r_0 = a_0 = \mu / R_0^2$. Consequently, one may compute
\begin{equation}
    R_0 = \left[ \mu  \left(\frac{\beta_j - \beta_i}{t_j - t_i} \right) \right]^{1/3}
\end{equation}
where the change in heading $\beta_j - \beta_i$ may be found directly from the heading directions $\bs'_i$ as 
\begin{equation}
    \beta_j - \beta_i = \text{atan2}\left[ \hat{\bw}^T (\bs'_i \times \bs'_j), \, (\bs'^T_i \bs'_j) \right]
\end{equation}
with $0 \leq \beta_j - \beta_i \leq 2 \pi$. The initial guess for $R_0$ may be computed from just a single heading pair, as an average of all the heading pairs, or from some other reasonable scheme.

Therefore, suppose the vector of parameters to estimate is given by $\bx$
\begin{equation}
    \label{eq:DefParamVector}
    \bx = \begin{bmatrix} x_1 \\ x_2 \\ x_3 \end{bmatrix} = \begin{bmatrix} R \\ c'_1 \\ c'_2 \end{bmatrix}
\end{equation}
Then a good initial guess for this vectors is
\begin{equation}
    \label{eq:HodographInitialGuess}
    \bx^{(0)} = \begin{bmatrix} R_0 \\ 0 \\ 0 \end{bmatrix}
\end{equation}
This initial guess has been found to reliably converge to the correct solution for even high eccentricity ($e \sim 0.9$) orbits.

\subsubsection{Developing the Iterative Scheme}
Beginning with an initial guess for the parameter vector $\bx^{(0)}$ from Eq.~\eqref{eq:HodographInitialGuess}, an iterative nonlinear least squares algorithm may be used to find the hodograph. It is most convenient to iterate on $\bx$ to minimize the differences between the measured and predicted time-of-flight. To make this explicit, first consider the equation for eccentric anomaly from Eq.~\eqref{eq:EccAnomalyHodographATAN2}. In the case where $\bc'$ is very small (as happens with the proposed initial guess), the orientation of the perifocal basis vectors becomes ill-defined and the expression for eccentric anomaly becomes atan2(0,0). Therefore, for the purposes of numeric computation, introduce a new vector $\bd'$ defined as
\begin{equation}
    \bd' = \left\{ \begin{array}{l l} 
                \bc' & e = \| \bc' \| / R > e_{min} \\
                \left[1,0,0\right]^T & \text{otherwise}
                \end{array}\right.
\end{equation}
where $e_{min}$ is the threshold where an orbit is considered effectively circular. This work uses $e_{min} \approx 0.001$, but the reader is free to make a different choice. The vector $\bd'$ is used to define a consistent reference direction in the equation for eccentric anomaly such that
\begin{equation}
    \label{eq:EccAnomalyHodographATAN2numeric}
    E_i = \text{atan2}\left[ \left(\sqrt{R^2-\bc'^T \bc'} \, \left[ \hat{\bw}^T ( \bd' \times \bs'_i) \right] \right) ,\, \left( R \, \bd'^T \bs'_i \right) \right] 
\end{equation}
Since real orbits are only approximately conics, the error introduced by using an arbitrary direction $\bd'=\left[1,0,0\right]^T$ for very nearly circular orbits is inconsequential in practice. 

With numerical issues for low eccentricity orbits resolved, observe that the mean anomaly $M_i$ may be computed in terms of the parameter $\bx$ (which contains $R$ and $\bc'$) and the heading $\bs'_i$ by substitution Eqs.~\eqref{eq:EccFromHodograph} and \eqref{eq:EccAnomalyHodographATAN2numeric} into Eq.~\eqref{eq:KeplersEqn}. Thus, one can write a function $g$ such that
\begin{equation}
\begin{aligned}
    M_i & =  g(\bx,\bs'_i)\\
    \; & = \text{atan2}\left[ \left(\sqrt{R^2-\bc'^T \bc'} \, \left[ \hat{\bw}^T ( \bd' \times \bs'_i) \right] \right) ,\, \left( R \, \bd'^T \bs'_i \right) \right] 
    \\
    & \quad + \frac{\| \bc' \|}{R} \sin \left\{ \text{atan2}\left[ \left(\sqrt{R^2-\bc'^T \bc'} \, \left[ \hat{\bw}^T ( \bd' \times \bs'_i) \right] \right) ,\, \left( R \, \bd'^T \bs'_i \right) \right] \right\}
\end{aligned}
\end{equation}
Moreover, one can construct the function $f$ describing the time of flight (for $t_j > t_i$) as
\begin{equation}
    \begin{aligned}
   \Delta t_{ij} = t_j - t_i & = \frac{1}{n}( 2 \pi k + M_j - M_i) \\
   & = \frac{\mu}{(R^2 -\bc'^T \bc')^{3/2}}\left[ 2 \pi k + g(\bx,\bs'_j) - g(\bx,\bs'_i) \right] \\
   & = f(\bx,\bs'_i,\bs'_j)
   \end{aligned}
\end{equation}
where $k$ is the number of times the orbiting body passes through periapsis between times $t_i$ and $t_j$.
Selecting $k$ is straightforward if all the heading observations occur within one orbital period (i.e., when $\Delta t_{ij} < P = 2 \pi / n$). In such a situation, set $k=0$ and compute $f(\bx,\bs'_i,\bs'_j)$. If $f<0$, set $k=1$ and recompute $f$. This check works since $\Delta t_{ij} > 0$ by construction. Schemes for selecting $k$ under more complicated situations are left to the reader.

The time-of-flight predicted by $f(\bx,\bs'_i,\bs'_j)$ may be compared to the measured time-of-flight to iteratively refine the estimate of $
\bx$ (which, again, contains the parameters $R$ and $\bc'$). Taking a Taylor series expansion of the time-of-flight $\Delta t_{ij}$ about an estimated hodograph at iteration $m$ yields
\begin{equation}
   \Delta t_{ij} = \Delta t^{(m)}_{ij} + \frac{\partial f_{ij}}{\partial \bx} \delta \bx^{(m)} + \mathcal{O}( \| \delta \bx \|^2 )
\end{equation}
which may be rearranged to find
\begin{equation}
   \Delta t_{ij} -\Delta t^{(m)}_{ij} \approx \frac{\partial f_{ij}}{\partial \bx} \delta \bx^{(m)}
\end{equation}
Assuming four observations, construct and stack the $\binom{4}{2} = 6$ combinations of times to form the linear system
\begin{equation}
    \begin{bmatrix}
   \Delta t_{12} -\Delta t^{(m)}_{12}\\
   \Delta t_{13} -\Delta t^{(m)}_{13}\\
   \vdots \\
   \Delta t_{34} -\Delta t^{(m)}_{34}\\
   \end{bmatrix}
   \approx 
   \begin{bmatrix}
   \frac{\partial f_{12}}{\partial \bx}\\
   \frac{\partial f_{13}}{\partial \bx}\\
   \vdots \\
   \frac{\partial f_{34}}{\partial \bx}\\
   \end{bmatrix}
    \delta \bx^{(m)}
\end{equation}
which may be iteratively solved using the Levenberg-Marquardt algorithm \cite{Levenberg:1944,Marquardt:1963}. This same concept may be extended without modification to accommodate $m$ heading observations within a $\binom{m}{2} \times 3$ linear system. At each iteration the estimate of $\bx$ is updated according to 
\begin{equation}
    \bx^{(k+1)} = \bx^{(m)} + \delta \bx^{(m)}
\end{equation}
Once converged to the optimal solution $\bx^{\ast}$, the IOD solution for the hodograph parameters may be directly recovered using Eq.~\eqref{eq:CenterInverseMapping} and Eq.~\eqref{eq:DefParamVector}
\begin{equation}
    \begin{aligned}
     &R = \bx^{\ast}_1 \\
     &\bc = c'_1 \hat{\ba} + c'_2 \hat{\bb}
          = \bx^{\ast}_2 \, \hat{\ba} +\bx^{\ast}_3 \, \hat{\bb}
    \end{aligned}
\end{equation}
The perifocal basis (only well-defined for non-circular orbits) may also be computed as from $\bc$ and $\hat{\bw}$
\begin{subequations}
\begin{equation}
    \hat{\bq} = \bc / \| \bc \| \\
\end{equation}
\begin{equation}
    \hat{\bp} = \hat{\bq} \times \hat{\bw} \\
\end{equation}
\end{subequations}
The conventional orbital elements may be computed from the hodograph parameters and perifocal basis vectors if desired.

\section{Numerical Example}
The efficacy of the proposed heading-only IOD scheme will now be discussed within the context of a simple example. Consider a spacecraft in low lunar orbit (LLO) with orbital elements given in Table~\ref{tab:OE}.
At four points around the orbit, a visual odometry algorithm (e.g., Ref.~\cite{Christian:2021b}) is used to obtain the heading measurements $\{ \bs_i \}_{i=1}^4$ given in Table~\ref{tab:Headings}. For context, Table~\ref{tab:RefAngles} includes the true anomaly, heading angle, eccentric anomaly, and time from periapsis at these four points. Given only the heading vectors from Table~\ref{tab:Headings}, the procedure from Section~\ref{Sec:IODMethod} is followed and results in the iteration history shown in Table~\ref{tab:IterHist}. The initial guess of $\bx^{(0)} = [1.4989,0,0]^T$ obtained using the scheme from Section~\ref{Sec:InitialGuess} was found to converge to the true solution in five iterations.

\begin{table}[h!]
    \centering
    \caption{Orbital elements of example spacecraft in low lunar orbit.}
    \label{tab:OE}
    \begin{tabular}{l c}
        \hline
        \hline
        Parameter & Value \\
        \hline
        Semi-major axis, $a$ & 2173.4 km \\
        Eccentricity, $e$ & 0.15 \\
        Inclination, $i$ & 65 deg\\
        Right ascension of ascending node, $\Omega$ & 70 deg \\
        Argument of periapsis, $\omega$ & 20 deg \\
        \hline
        \hline
    \end{tabular}
\end{table}

\begin{table}[h!]
    \centering
    \caption{Heading observations (as unit vectors) that are used to solve the example IOD problem.}
    \label{tab:Headings}
    \begin{tabular}{l c c c}
        \hline
        \hline
        Heading & \multicolumn{3}{c}{Unit Vector} \\
        \cline{2-4}
        Observation & $x$ & $y$ & $z$ \\
        \hline
        $\bs_1$ & -0.5028 &  -0.2557  &  0.8257 \\
        $\bs_2$ & -0.3918 &  -0.9122 &   0.1204 \\
        $\bs_3$ & 0.2052  & -0.5448 &  -0.8131 \\
        $\bs_4$ & 0.3900  &  0.9135 &  -0.1158 \\
        \hline
        \hline
    \end{tabular}
\end{table}

\begin{table}[h!]
    \centering
    \caption{Reference details for each heading observation from Table~\ref{tab:Headings}. This information is for context and is not known to the IOD algorithm.}
    \label{tab:RefAngles}
    \begin{tabular}{l c c c c}
        \hline
        \hline
        Heading & True Anomaly & Heading Angle & Eccentric Anomaly & Time from Periapsis\\
        Observation & $\theta_i$ [deg] & $\beta_i$ [deg] & $E_i$ [deg] & $t_i-t_0$ [min] \\
        \hline
        $\bs_1$ &  5.00  &  4.35  &  4.30& 1.54 \\
        $\bs_2$ & 70.00  & 62.36  & 62.09 & 22.94 \\
        $\bs_3$ & 140.00 & 133.78 & 134.11 & 53.85 \\
        $\bs_4$ & 235.00 & 242.66 & 242.39 & 105.24 \\
        \hline
        \hline
    \end{tabular}
\end{table}

\begin{table}[h!]
    \centering
    \caption{Iteration history for heading-based IOD method using only the four perfect (noise-free) heading observations from Table~\ref{tab:Headings}.}
    \label{tab:IterHist}
    \begin{tabular}{l c c c c}
        \hline
        \hline
        Iteration & & \multicolumn{3}{c}{Parameter Vector $\bx^{(m)}$} \\
        \cline{3-5}
        Number, $m$ & Residual& $R$ & $c'_1$ & $c'_2$ \\
        \hline
           0    &  314791 & 1.4989   &      0.0   &      0.0     \\
            1   &  55148 & 1.4752  & -0.0937 &  -0.0937   \\
            2   &  21355 & 1.5067  & -0.2110  &  0.0071    \\
            3   &  10.38  & 1.5188 &  -0.2267  &  0.0170    \\
            4   &  $4.53 \times 10^{-6}$ & 1.5191 &  -0.2272  &  0.0173 \\  
            5   &  $6.68 \times 10^{-19}$ & 1.5191 &  -0.2272 &   0.0173\\
        \hline
        \hline
    \end{tabular}
\end{table}

To convert the results for $\bc'$ from Table~\ref{tab:IterHist} into the final hodograph center $\bc$, observe that
$$
    \hat{\ba}  = \begin{bmatrix}
    0.5028 \\
    0.2557 \\
   -0.8257
    \end{bmatrix}
    \quad \text{and} \quad
    \hat{\bb}  = \begin{bmatrix}
    0.1479 \\
    0.9157 \\
    0.3737
    \end{bmatrix}
$$
such that, from Eq.~\eqref{eq:CenterInverseMapping}, one finds
$$
   \bc = c'_1 \hat{\ba} + c'_2 \hat{\bb} = 
   \begin{bmatrix}
    -0.1117 \\
   -0.0423 \\
    0.1941
    \end{bmatrix}
$$
which the reader can verify is the correct vector $\bc$ that is computed from the orbital elements in Table~\ref{tab:OE}. Without heading errors, this algorithm converges to machine precision (or to the tolerance chosen for the Levenberg-Marquardt algorithm).

There are a number of factors that make the perfect (noise-free) situation described above uncommon in practice. First, heading measurements from image-based visual odometry are not perfect. Empirical evidence (e.g, from Ref.~\cite{Christian:2021b}) suggests that heading errors on the order of about 0.1--1.0 deg should be expected. Second, the central body (Moon in this example) is rotating and surface features are not inertially fixed. The siderial rotational period of the Moon is about 27.3 days (or about $2.36 \times 10^6$ seconds) leading to an equatorial velocity of about 4.6 m/s. The apoapsis velocity of this example orbit is 1.3 km/s (slowest point), meaning the maximum heading error introduced by neglecting the Moon's rotation is about 0.2 deg (and usually much less). 

\begin{figure}[b!]
    	\centering
    	\includegraphics[width=0.95 \textwidth]{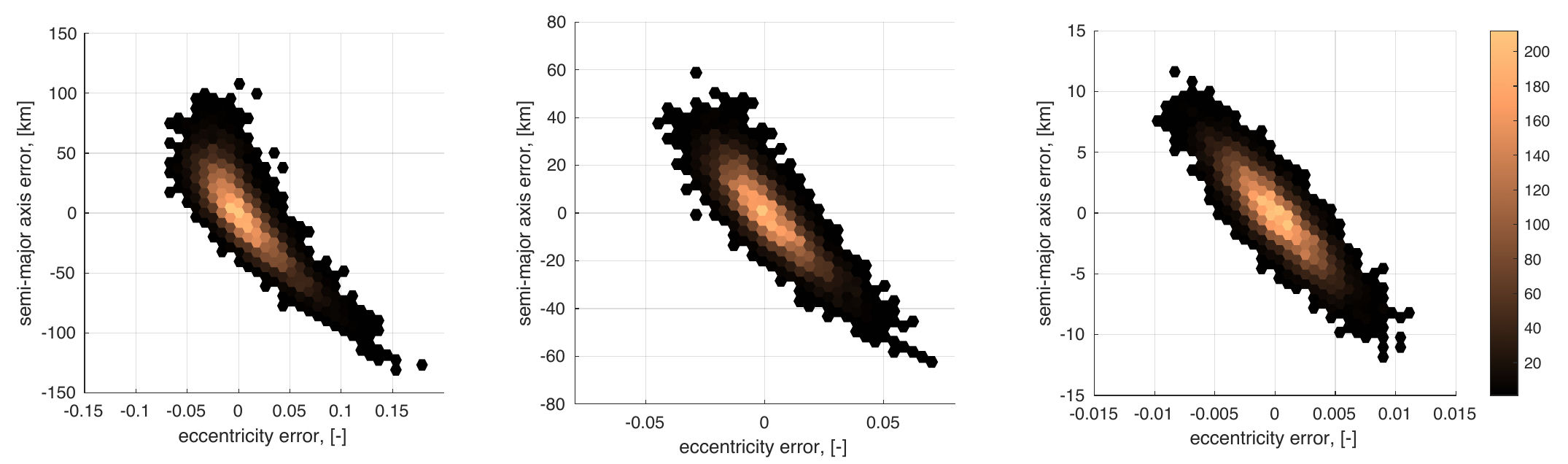}
    	\caption{IOD errors from 10,000-run Monte Carlo for heading errors of 1 deg (left), 0.5 deg (middle), and 0.1 deg (right) with the four heading observations from Table~\ref{tab:Headings}}.
    	\label{fig:FourObsErrors}
\end{figure}

\begin{figure}[b!]
    	\centering
    	\includegraphics[width=0.95 \textwidth]{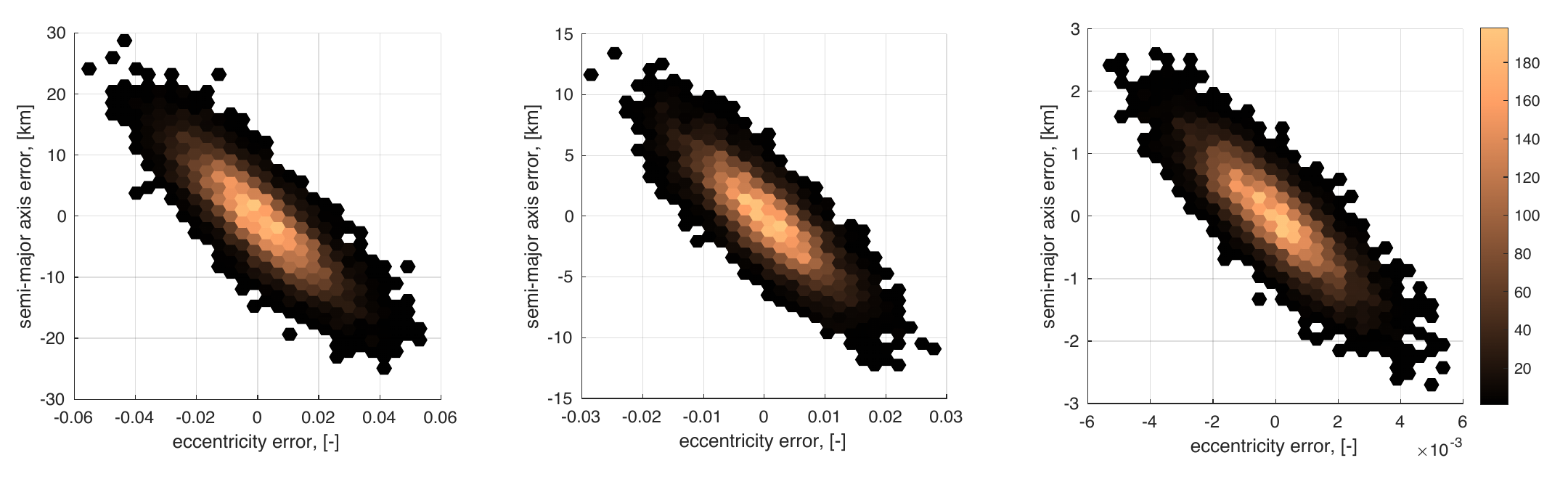}
    	\caption{IOD errors from 10,000-run Monte Carlo for heading errors of 1 deg (left), 0.5 deg (middle), and 0.1 deg (right) with ten evenly spaced heading observations between $\theta_1 = 15$ deg and $\theta_{10} = 330$ deg. }
    	\label{fig:TenObsErrors}
\end{figure}

Therefore, consider the IOD performance of the four-heading example above with three noise levels on the measured heading directions: 1.0 deg, 0.5 deg, and 0.1 deg. A 10,000-run Monte Carlo analysis was performed for each scenario, and plots of the resulting semi-major axis error and eccentricity error were computed. These are shown in Fig.~\ref{fig:FourObsErrors}, where IOD performance improves as measurement noise decreases (left to right). As a comparison, suppose ten evenly-spaced (in true anomaly) headings were collected between $\theta_1 = 15$ deg and $\theta_{10} = 330$ deg. The exact same IOD algorithm may be performed (the method of this paper is not limited to just four headings) and the results are shown in Fig.~\ref{fig:TenObsErrors}. The results of both scenarios (four heading measurements and ten heading measurements) are summarized and compared in Table~\ref{tab:PerfMC}. As expected, more observations leads to better IOD performance and increasing measurement noise leads to worse IOD performance.

\begin{table}[t!]
    \centering
    \caption{Summary of heading-only IOD performance for 10,000-run Monte Carlo simulations shown in Fig.~\ref{fig:FourObsErrors} and Fig.~\ref{fig:TenObsErrors}. Error statistics represent one standard deviation ($1\sigma$) values.}
    \label{tab:PerfMC}
    \begin{tabular}{l c c c }
        \hline
        \hline
        Number of & Heading Error  & $a$ Error & $e$ Error\\
        Observations & [deg] & [km]  & [-]\\
        \hline
        4  & 1.0 & 31.2721 & 0.0287 \\
        4   & 0.5 & 15.4026 & 0.0140 \\
        4    &0.1 & 3.0635 & 0.0027\\
        \hline
        10   & 1.0 & 7.1623 & 0.0145\\
        10 &    0.5 & 3.5655 & 0.0072 \\
        10   &    0.1 & 0.7174 & 0.0015 \\
        \hline
        \hline
    \end{tabular}
\end{table}

\section{Conclusions}
This work introduces a new class of initial orbit determination (IOD) problem. In this new problem, one seeks to find the orbit of a body using only four (or more) observations of heading direction. Such heading observations are produced by camera-based navigation systems that operate on the principle of visual odometry. A simple iterative solution to this problem is constructed using the orbital hodograph and is found to quickly converge to the true solution in the absence of noise. As one would expect, it is found that IOD performance degrades with increasing measurement noise and improves with an increasing number of measurements. The author hopes that others find this new class of IOD problem interesting and will contribute their own solutions to the literature.

\section*{Funding Sources}
This work was partially supported by the Air Force Office of Scientific Research (AFOSR) under the Space University Research Initiative (SURI) grant FA9550-22-1-0092 (grant principal investigator: J. Crassidis from University at Buffalo, The State University of New York).

\bibliography{sample}

\end{document}